\documentclass[11pt,a4paper]{amsart}
\usepackage[utf8]{inputenc}
\usepackage{amsmath}
\usepackage{amsfonts}
\usepackage{amssymb}
\usepackage{graphicx}
\usepackage{booktabs}

\usepackage[
left=2.0cm, 
right=2.0cm, 
top=1.8cm, 
bottom=1.8cm 
]{geometry}

\newcommand{\D}{{\mathop{}~\!\mathrm{d}}} 
\newcommand{\R}{\mathbb{R}}
\newcommand{\Q}{\mathbb{Q}}

\newcommand{\N}{\mathbb{N}}

\newcommand{\PP}{\mathbb{P}}

\numberwithin{equation}{section}  
\newtheorem{defn}{Definition}[section]

\newtheorem{rem}[defn]{Remark}
\newtheorem{thm}[defn]{Theorem}
\newtheorem{prop}[defn]{Proposition}

\newtheorem{lem}[defn]{Lemma}

\usepackage{hyperref}
\usepackage{xcolor}
\hypersetup{
    colorlinks,
    linkcolor={red!50!black},
    citecolor={blue!50!black},
    urlcolor={blue!80!black}
}

\title[On the stability of the martingale optimal transport problem]{On the stability of the martingale optimal transport problem: A set-valued map approach}

\author{Ariel Neufeld \and Julian Sester}
\address{NTU Singapore, Division of Mathematical Sciences; 21 Nanyang Link, Singapore 637371}
\begin{document}
\vspace{-0.5cm}

\begin{abstract}
Continuity of the value of the martingale optimal transport problem on the real line w.r.t.\ its marginals was recently established in \cite{backhoff2019stability} and \cite{wiesel2019continuity}. We present a new perspective of this result using the theory of set-valued maps. In particular, using results from \cite{beiglbock2021approximation}, we show that the set of martingale measures with fixed marginals is continuous, i.e., lower- and upper hemicontinuous, w.r.t.\ its marginals. Moreover, we establish compactness of the set of optimizers as well as upper hemicontinuity of the optimizers w.r.t.\ the marginals. 
\end{abstract}

\begin{minipage}[t]{\textwidth}
\vspace*{-4\baselineskip}
\maketitle
\end{minipage}

\section{Introduction}
\noindent
The martingale optimal transport problem (as introduced in \cite{beiglbock2013model}) consists, given a measurable function $\Phi:\R^2 \rightarrow \R$, in solving
\begin{equation}\label{eq_mot_mu1mu2}
m(\mu_1,\mu_2):=\sup_{\Q \in \mathcal{M}(\mu_1,\mu_2)}\int_{\R^2}\Phi(x_1,x_2)\D \Q(x_1,x_2),
\end{equation}
where $\mathcal{M}(\mu_1,\mu_2)$ describes the set of martingale measures on $\R^2$ with fixed marginals $\mu_1$ and $\mu_2$.

Recently, the stability of the martingale optimal transport problem w.r.t.\ its marginal distributions, i.e., the question whether the solutions $m(\mu_1,\mu_2)$  and $m(\widetilde{\mu_1},\widetilde{\mu_2})$ of martingale optimal transport problems are close, whenever the marginals $\mu_i$, and $\widetilde{\mu_i}$, $i=1,2$, are close in the Wasserstein-distance, was studied in several research papers. One particular reason for the importance of a positive answer to this question is that in this case martingale transport problems involving finitely supported marginal distribution that are close to the original marginals yield solutions that are close to the original value. Since these approximated solutions can be computed with tractable linear programming methods (compare e.g. \cite{guo2019computational} and \cite{henry2013automated}), the stability result builds an important theoretical foundation for the numerics of martingale optimal transport.

The stability result was indeed established for a Lipschitz-continuous cost function and a relaxed formulation of the transport problem in \cite{guo2019computational}, for particular cost functions and special marginals in \cite{juillet2016stability}, and eventually in a great generality by \cite{backhoff2019stability} and \cite{wiesel2019continuity}. Recently, \cite{beiglbock2021approximation} established a result which allows to approximate martingale measures with fixed marginals $(\mu_1,\mu_2)$ in the adapted Wasserstein-distance (compare \cite{backhoff2019stability}) by a sequence of approximating martingale measures with fixed marginals $(\mu_1^{(n)},\mu_2^{(n)})$ which converge in Wasserstein-distance to $(\mu_1,\mu_2)$ for $n \rightarrow \infty$.

In this work, we regard the problem from a new perspective by studying properties of set-valued maps. More precisely, we show that the set-valued map
\begin{equation}\label{eq_hemicontinuity_intro}
(\mu_1,\mu_2) \mapsto \mathcal{M}(\mu_1,\mu_2),
\end{equation}
is continuous, i.e., is upper hemicontinuous and lower hemicontinuous.
This in turn implies through an application of Berge's maximum theorem \cite{berge} the stability of the martingale optimal transport problem w.r.t.\ its marginals.
Moreover, we obtain compactness of the set of optimizers, i.e., of the set of two-dimensional martingale measures maximizing $m(\mu_1,\mu_2)$ as well as the upper-hemicontinuity of the correspondence mapping from marginals to these optimizers.

\section{Main Result}
We first introduce some notation. For every $n\in \N$ let $\mathcal{P}(\R^n)$ denote the set of all probability measures on $\R^n$. Then, we define the $1$-\emph{Wasserstein space}
\[
\mathcal{P}_1(\R^n):=\left\{ \PP \in \mathcal{P}(\R^n)~\middle|~\int_{\R^n} \textstyle \sum \limits_{i=1}^n|x_i|\D \PP(x_1,\dots,x_n) < \infty\right\}
\]
describing the probability measures with existing first moment. Moreover, we let
$$
C_{\operatorname{lin}}\left(\R^n\right):= \left\{f\in C\left(\R^n\right)~\middle|~ \sup_{(x_1,\dots,x_n)\in \R^n} \frac{|f(x_1,\dots,x_n)|}{1+\sum_{i=1}^n|x_i|}  < \infty\right\}
$$
denote the set of all continuous functions from $\R^n$ to $\R$ with linear growth.
Next, we define for every $n \in \N$ the set $\Pi(\mu_1,\mu_2)$ of all probability measures on $\R^{2n}$ with fixed marginal  distributions $\mu_1, \mu_2 \in \mathcal{P}_1(\R^n)$, also called couplings of $\mu_1,\mu_2$. 
We further define for $\PP,\Q \in \mathcal{P}_1(\R^n)$ the $1$-Wasserstein distance
\[
\mathcal{W}(\PP,\Q):=\inf_{\pi\in \Pi(\PP,\Q)}\int_{\R^n\times \R^n} \textstyle \sum \limits_{i=1}^n|x_i-y_i|\D\pi(x_1,\dots,x_n,y_1,\dots,y_n),
\]
as well as the sum of Wasserstein distances between measures $\PP_i,\Q_i \in \mathcal{P}_1(\R^n)$ for $i=1,\dots,m$
\[
\mathcal{W}^{\oplus}\big((\PP_1,\dots,\PP_m),(\Q_1,\dots,\Q_m)\big):=\sum_{i=1}^m \mathcal{W}(\PP_i,\Q_i).
\]

The set of all martingale measures on $\R^{2}$ with fixed marginal distributions $\mu_1, \mu_2 \in \mathcal{P}_1(\R)$ is then given by
\begin{equation}\label{eq_martingale_def}
\begin{aligned}
\mathcal{M}(\mu_1,\mu_2):=\bigg\{\Q \in \Pi(\mu_1,\mu_2)~\bigg|&~\int_{\R^2} H(x_1)(x_2-x_1)\D\Q(x_1,x_2)=0 \text{ for all } H \in C_b(\R)\bigg\},
\end{aligned}
\end{equation}
where $C_b(\R)$ denotes the set of bounded continuous functions from $\R$ to $\R$.
Moreover, we denote by $\mu_1 \preceq \mu_2$ the convex order of measures $\mu_1,\mu_2 \in \mathcal{P}_1(\R)$, which is characterized by the inequality $\int f \D \mu_1 \leq \int f \D \mu_2$ required to hold for all convex functions $f:\R\rightarrow \R$. Strassen's result (\cite[Theorem 8]{strassen1965existence}) implies\footnote{In \cite[Theorem 8]{strassen1965existence}, Strassen proves the existence of a martingale with prescribed marginal distributions.} that $\mathcal{M}(\mu_1,\mu_2) \neq \emptyset$ if and only if $\mu_1 \preceq \mu_2$.

When speaking of a \emph{correspondence} from a set $X$ to a set $Y$, denoted by $\varphi:X\twoheadrightarrow Y$, we refer to a set-valued function, i.e., a map that assigns each $x \in X$ a set $\varphi(x)\subseteq Y$, see also \cite[Chapter 17]{Aliprantis}. What continuity means for a correspondence is clarified in the following definition.
\begin{defn}\label{def_hemi}
Let $\varphi:X \twoheadrightarrow Y$ be a correspondence between two topological spaces.
\begin{itemize}
\item[(i)]
$\varphi$ is called \emph{upper hemicontinuous}, if $\{x \in X~|~\varphi(x) \subseteq A\}$ is open for all open sets $A \subseteq Y$.
\item[(ii)] $\varphi$ is called \emph{lower hemicontinuous}, if 
$\{x \in X~|~ \varphi(x) \cap A \neq \emptyset\}$ is open for all open sets $A \subseteq Y$.
\item[(iii)] We say $\varphi$ is continuous, if $\varphi$ is upper and lower hemicontinuous.
\end{itemize}
\end{defn}
Our main result is the following.
\begin{prop}\label{prop_main_result}
Let
\[
X:=\left\{(\mu_1,\mu_2)\in \mathcal{P}_1(\R)\times\mathcal{P}_1(\R) ~\middle|~\mu_1\preceq \mu_2\right\}
\]
be equipped with the topology induced by the convergence  w.r.t.\ $\mathcal{W}^{\oplus}$, and let 
\[
Y:=\mathcal{P}_1(\R^2)
\]
be equipped with the topology induced by the convergence w.r.t.\ $\mathcal{W}$.
Then, the correspondence
\begin{align*}
\varphi: X &\twoheadrightarrow Y\\
(\mu_1,\mu_2) &\mapsto \mathcal{M}(\mu_1,\mu_2)
\end{align*}
has non-empty, convex and compact images, and is continuous. As a consequence one obtains the following.
\begin{itemize}
\item[(i)]
For any $\Phi \in C_{\operatorname{lin}}\left(\R^2\right)$ the functional
\begin{align*}
m: X &\rightarrow \R \\
(\mu_1,\mu_2) &\mapsto \max_{\Q \in \mathcal{M}(\mu_1,\mu_2)} \int_{\R^2} \Phi(x_1,x_2) \D \Q(x_1,x_2)
\end{align*}
is continuous.
\item[(ii)] For any $(\mu_1,\mu_2) \in X$ and $\Phi \in C_{\operatorname{lin}}\left(\R^2\right)$ the set of optimizers 
\[
\mathcal{Q}^*(\mu_1,\mu_2):=\left\{\Q^* \in \mathcal{M}(\mu_1,\mu_2)~\middle|~\max_{\Q \in \mathcal{M}(\mu_1,\mu_2)}\int_{\R^2} \Phi(x_1,x_2) \D \Q(x_1,x_2) = \int_{\R^2} \Phi(x_1,x_2) \D \Q^*(x_1,x_2)\right\}
\]
is non-empty and compact.
\item[(iii)]
For any $\Phi \in C_{\operatorname{lin}}\left(\R^2\right)$ the correspondence
\begin{align*}
\mu:X &\twoheadrightarrow Y \\
(\mu_1,\mu_2) &\mapsto \mathcal{Q}^*(\mu_1,\mu_2)
\end{align*}
is upper hemicontinuous.
\end{itemize}
\end{prop}

\begin{rem}
\begin{itemize}

\item[(i)]
Let $X,Y$ denote topological spaces.
If a correspondence $X \supseteq \Theta \ni \theta \mapsto \mathcal{S}(\theta) \subseteq Y$ is single-valued, i.e., $\mathcal{S}(\theta)=\{s_\theta\}$ for all $\theta \in \Theta$, then upper hemicontinuity of the set-valued map $\Theta \ni \theta \mapsto \mathcal{S}(\theta) $ is equivalent to the continuity of the (single-valued) map $X \supseteq \Theta \ni \theta\mapsto s_\theta \in Y$ as a function from topological space $\Theta\subseteq X$ to $Y$, compare e.g. \cite[Lemma 17.6]{Aliprantis}. Thus, if the set of optimizers  $\mathcal{Q}^*(\mu_1,\mu_2)=\{\Q^*(\mu_1,\mu_2)\}$ is a singleton for all $(\mu_1,\mu_2) \in X$, we have by Proposition~\ref{prop_main_result} that the map
\[
\left\{(\mu_1,\mu_2)\in \mathcal{P}_1(\R)\times\mathcal{P}_1(\R) ~\middle|~\mu_1\preceq \mu_2\right\}\ni (\mu_1,\mu_2) \mapsto \Q^*(\mu_1,\mu_2) \in \mathcal{P}_1(\R^2)
\]
is continuous, if $\Phi \in C_{\operatorname{lin}}(\R^2)$. It is possible to state explicit conditions on the cost function $\Phi$ which imply uniqueness of the associated optimizers. These conditions were among others studied in detail in  \cite{beiglbock2016problem}, \cite{henry2016explicit}, \cite{hobson2015robust},  and \cite{hobson2019left}.
\item[(ii)]
If a correspondence is non-empty, convex, compact, and continuous, then there exist results enabling a (numerically tractable) approximation of the correspondence, whenever its image is contained in a finite-dimensional normed vector space, compare for example \cite{dudov2007approximation} and \cite{dynfarkhi2006}. The correspondence $(\mu_1,\mu_2) \mapsto \mathcal{M}(\mu_1,\mu_2)$ fulfils, according to Proposition~\ref{prop_main_result}, all of the above mentioned requirements, but its image is not contained in a finite-dimensional normed vector space.  
If a suitable approximation would be possible also for this measure-valued correspondence, then such an approximation could allow to approximate the martingale optimal transport problem \eqref{eq_mot_mu1mu2} under consideration.
Moreover, the continuity property of the set valued map $(\mu_1,\mu_2) \mapsto \mathcal{M}(\mu_1,\mu_2)$ could be fruitful to analyze (lower hemi-) continuity of the set valued map of optimizers w.r.t.\ the given marginals, i.e. $(\mu_1,\mu_2) \mapsto  \mathcal{Q}^*(\mu_1,\mu_2)$, compare, e.g., with \cite{ kien2005lower}, which in turn leads to the existence of continuous selectors, see, e.g., \cite{michael1956continuous}.
\\
We leave these questions open for future research.
\item[(iii)]
The martingale optimal transport problem introduced in \cite{beiglbock2013model} is formulated as a minimization problem. We decided however to use the maximization formulation, since for financial applications it is often of higher interest to compute upper price bounds, compare also the formulations in \cite{beiglbock2017monotone}, \cite{guo2019computational}, \cite{henry2013automated}, \cite{henry2016explicit} which all consider maximization problems. By the relation 
\[
\inf_{\Q \in \mathcal{M}(\mu_1,\mu_2)}\int\Phi(x_1,x_2) \D \Q(x_1,x_2) = -\sup_{\Q \in \mathcal{M}(\mu_1,\mu_2)}\int -\Phi(x_1,x_2) \D \Q(x_1,x_2),
\] the corresponding results from Proposition~\ref{prop_main_result} also remain valid for the minimization problem.
\item[(iv)]
The non-emptiness of $\mathcal{Q}^*(\mu_1,\mu_2)$ can be established when only requiring upper semicontinuity of $\Phi$, see e.g. \cite[Theorem 1.1.]{beiglbock2013model}. However, for the continuity of the functional ${m}$ it seems to be necessary to require that $\Phi$ is continuous, compare also the main results from \cite{backhoff2019stability} and \cite{wiesel2019continuity} which are formulated with the same requirement on $\Phi$.
\item[(v)] Recently, Br\"{u}ckerhoff and Juillet provided in \cite{bruckerhoff2021instability} a counterexample implying that Proposition~\ref{prop_main_result} is indeed only valid for one-dimensional marginal distributions $\mu_1,\mu_2 \in \mathcal{P}_1(\R)$, but not for marginal distributions from $\mathcal{P}_1(\R^d)$ with $d \geq 2$.
\end{itemize}
\end{rem}

\section{Proof of Proposition~\ref{prop_main_result}}
\subsection{Preliminaries and useful Results}
For the proof of the main result we make use of the following results.
The first result is Berge's maximum theorem. Note that for a correspondence $\varphi:X \twoheadrightarrow Y$ the graph is defined as 
\[
\operatorname{Gr}\varphi: = \left\{(x,y) \in X \times Y ~\middle|~y \in \varphi(x)\right\}.
\]
\begin{thm}[\cite{Aliprantis}, Theorem 17.31]\label{thm_berge}
Let $\varphi:X \twoheadrightarrow Y$ be a continuous correspondence between topological spaces with non-empty compact values, and suppose that $f: \operatorname{Gr}\varphi \rightarrow \R$ is continuous. Then the following holds.
\begin{itemize}
\item[(i)]
The function 
\begin{align*}
m: X &\rightarrow \R\\
x&\mapsto\max_{y \in \varphi(x)}f(x,y)
\end{align*}
is continuous.
\item[(ii)]
The correspondence
\begin{align*}
c: X &\twoheadrightarrow Y \\
x &\mapsto \left\{y \in \varphi(x)~\middle|~f(x,y)=m(x)\right\}
\end{align*}
has non-empty, compact values.
\item[(iii)]
If $Y$ is Hausdorff, then $c$ is upper hemicontinuous.
\end{itemize}
\end{thm}
By Definition~\ref{def_hemi} a correspondence is called continuous if it is upper hemicontinuous and lower hemicontinuous. These notions are characterized through the following two lemmas.

\begin{lem}[\cite{Aliprantis}, Theorem 17.20]\label{lem_upper_hemi}
Assume that the topological space $X$ is first countable and that $Y$ is metrizable. Then, for a correspondence $\varphi:X \twoheadrightarrow  Y$  the following statements are equivalent.
\begin{itemize}
\item[(i)] The correspondence $\varphi$ is upper hemicontinuous and $\varphi(x)$ is compact for all $x\in X$.
\item[(ii)]
For any $x\in X$, if a sequence $\left((x^{(n)},y^{(n)})\right)_{n \in \N} \subseteq \operatorname{Gr}\varphi$ satisfies $x^{(n)} \rightarrow x$ for $n \rightarrow \infty$, then there exists a subsequence $\left(y^{(n_k)}\right)_{k \in \N}$ with $y^{(n_k)} \rightarrow y  \in \varphi(x)$ for $k \rightarrow \infty$.
\end{itemize}
\end{lem}

\begin{lem}[\cite{Aliprantis}, Theorem 17.21]\label{lem_lower_hemi}
For a correspondence $\varphi:X \twoheadrightarrow  Y$ between first countable topological spaces the following statements are equivalent.
\begin{itemize}
\item[(i)]
The correspondence $\varphi$ is lower hemicontinuous.
\item[(ii)]
For any $x\in X$, if $x^{(n)} \rightarrow x$ for $n \rightarrow \infty$, then for each $y \in \varphi(x)$ there exists a subsequence $\left(x^{(n_k)}\right)_{k \in \N}$ and elements $y^{(k)} \in \varphi\left(x^{(n_k)}\right)$ for each $k\in \N$ such that $y^{(k)} \rightarrow y$ for $k \rightarrow \infty$.
\end{itemize}
\end{lem}

The following result from \cite{beiglbock2021approximation} turns out to be crucial for our arguments. Note that the adapted $1$-Wasserstein distance between two measures $\Q,\Q'$ (with first marginals $\mu_1$ and $\mu_1'$, respectively) is defined as 
\[
\mathcal{AW}(\Q,\Q'):= \inf_{\pi \in \Pi(\mu_1,\mu_1')}\int_{\R^2} |x_1-x_1'| +\mathcal{W}(\Q_{x_1},\Q'_{x_1'})\D \pi(x_1,x_1'),
\]
where $\Q_{x_1},\Q'_{x_1'}$ denote the disintegration of the measures $\Q$ and $\Q'$ respectively w.r.t.\ their first marginals, i.e., 
$\Q = \mu_1 \otimes \Q_{x_1}$ and $\Q' = \mu_1' \otimes \Q'_{x_1'}$.
\begin{thm}[\cite{beiglbock2021approximation}, Theorem 2.5]\label{thm_beigl_continuity}
Let $\mu_1, \mu_2\in \mathcal{P}_1(\R)$, and consider $\mu_1^{(k)}, \mu_2^{(k)}\in \mathcal{P}_1(\R)$, $\mu_1^{(k)} \preceq \mu_2^{(k)}$,  $k \in \N$ with  
\[
\mathcal{W}^{\oplus}\left((\mu_1^{(k)}, \mu_2^{(k)}),(\mu_1,\mu_2)\right) \rightarrow 0 \text{ for }k \rightarrow \infty.
\]
Let $\Q \in \mathcal{M}(\mu_1,\mu_2)$. Then there exists a sequence $(\Q^{(k)})_{k \in \N}$ with $\Q^{(k)} \in \mathcal{M}(\mu_1^{(k)},\mu_2^{(k)})$ for all $k \in \N$ s.t.
\[
\mathcal{AW}(\Q^{(k)},\Q)\rightarrow 0 \text{ for } k \rightarrow \infty.
\]
\end{thm}
\subsection{Proof of Proposition~\ref{prop_main_result}}
\begin{proof} Let $\Phi \in C_{\operatorname{lin}}\left(\R^2\right)$ and define on 
\[
\operatorname{Gr}\varphi: =\big\{\left(\mu_1,\mu_2,\Q\right)~\big|~\mu_1,\mu_2 \in \mathcal{P}_1(\R),\mu_1\preceq \mu_2, \Q \in \mathcal{M}(\mu_1,\mu_2)\big\}
\]
the function
\begin{align*}
f: \operatorname{Gr}\varphi &\rightarrow \R \\
  \left(\mu_1,\mu_2,\Q\right) &\mapsto \int_{\R^2} \Phi(x_1,x_2) \D \Q(x_1,x_2).
\end{align*}
We aim at showing
\begin{itemize}
\item[1.)] $\varphi$ is non-empty valued and convex. 
\item[2.)] $\varphi$ is upper hemicontinuous and compact valued.
\item[3.)] $\varphi$ is lower hemicontinuous.
\item[4.)] $f$ is continuous.
\end{itemize}
If the above listed requirements are established, then the result  of Proposition~\ref{prop_main_result} follows directly by Berge's maximum theorem stated in Theorem~\ref{thm_berge}. Before proving the requirements we note that, according to  \cite[Definition 6.8~(iv) and Theorem 6.9]{villani2008optimal}, for $n\in \N$,
the convergence $\lim_{k\rightarrow\infty}\mathcal{W}(\PP,\PP^{(k)})= 0$  of a sequence $\left(\PP^{(k)}\right)_{k \in \N} \subseteq \mathcal{P}_1(\R^n)$ to some limit $\PP \in \mathcal{P}_1(\R^n)$ is equivalent to
\begin{equation}\label{eq_wassersteinspace_convergence}
\lim_{k\rightarrow\infty}\int_{\R^n} f(x) \D \PP^{(k)}(x) =\int_{\R^n} f(x) \D \PP(x) \text{ for all } f \in C_{\operatorname{lin}}(\R^n).
\end{equation}
Moreover, note that with the assigned metrics  $\mathcal{W}^{\oplus}$  and $\mathcal{W}$ respectively, the spaces $X$ and $Y$ are indeed first countable, thus Lemma~\ref{lem_upper_hemi} as well as Lemma~\ref{lem_lower_hemi} are applicable. Further, $Y$ is Hausdorff (since $\mathcal{W}$ is a metric on $\mathcal{P}_1(\R^2)$) as required in Theorem~\ref{thm_berge}.
\begin{itemize}
\item[1.)] Pick some $(\mu_1,\mu_2) \in X$. Then we have by definition of $X$ that $\mu_1 \preceq \mu_2$, which implies due to the well-known result from \cite[Theorem 8]{strassen1965existence} that $\varphi\left((\mu_1,\mu_2)\right)=\mathcal{M}(\mu_1,\mu_2) \neq \emptyset$. Moreover, note that the convexity of $\mathcal{M}(\mu_1,\mu_2)$ follows by definition.
%
\item[2.)]
To show the upper hemicontinuity of $\varphi$ and the compactness of the image we apply Lemma~\ref{lem_upper_hemi}. Thus, let $(\mu_1,\mu_2) \in X$ and consider a sequence $\left(\mu_1^{(n)},\mu_2^{(n)},\Q^{(n)}\right)_{n \in \N} \subseteq  \operatorname{Gr}(\varphi)$ such that 
\begin{equation}\label{eq_wasserstein_marginal_convergence1}
\lim_{n \rightarrow \infty}\mathcal{W}^{\oplus}\left((\mu_1^{(n)},\mu_2^{(n)}),(\mu_1,\mu_2)\right)= 0.
\end{equation}
Observe that by Prokhorov's theorem and by the weak convergence implied through \eqref{eq_wasserstein_marginal_convergence1} the sets $\{\mu_1^{(n)},n \in \N\}$ and $\{\mu_2^{(n)},n \in \N\}$ are tight. Denote by 
\[
\widetilde{\Pi}:=\Pi\left(\{\mu_1^{(n)},n \in \N\},\{\mu_2^{(n)},n \in \N\}\right) \subseteq  \mathcal{P}(\R^2)
\]the set of probability measures on $\R^2$ with first marginal in $\{\mu_1^{(n)},n \in \N\}$ and second marginal in $\{\mu_2^{(n)},n \in \N\}$. The set $\widetilde{\Pi}$ is tight according to \cite[Lemma 4.4]{villani2008optimal}. 
Thus, since $(\Q^{(n)})_{n\in \N} \subseteq  \widetilde{\Pi}$, according to Prokhorov's theorem,  there exists a subsequence $(\Q^{(n_k)})_{k\in \N}$ which converges weakly to some $\Q \in \mathcal{P}(\R^2)$.
We claim that $\Q\in \Pi(\mu_1,\mu_2)$. Indeed, let $\pi_i:\R^2 \rightarrow \R$ denote the projection on the $i$-th component for $i=1,2$. Then, by the continuous mapping theorem, it holds for $i=1,2$ that
\[
\Q^{(n_k)}\circ \pi_i^{-1} \rightarrow \Q\circ \pi_i^{-1} \text{ weakly for } k \rightarrow \infty.
\]
Additionally, by \eqref{eq_wasserstein_marginal_convergence1} and \eqref{eq_wassersteinspace_convergence}, it holds
\[
\Q^{(n_k)}\circ \pi_i^{-1} =\mu_i^{(n_k)} \rightarrow \mu_i \text{ weakly for } k \rightarrow \infty.
\]
Thus, $\Q \circ \pi_i^{-1}=\mu_i$ for $i=1,2$, i.e., $\Q \in \Pi(\mu_1,\mu_2)$.
Moreover, \eqref{eq_wasserstein_marginal_convergence1} together with \eqref{eq_wassersteinspace_convergence} implies 
\begin{align*}
\lim_{k\rightarrow \infty} \int_{\R^2}|x_1|+|x_2|\D\Q^{(n_k)}(x_1,x_2) &=\lim_{k\rightarrow \infty} \bigg(\int_{\R} |x_1| \D\mu_1^{(n_k)}(x_1)+\int_{\R} |x_2| \D\mu_2^{(n_k)}(x_2)\bigg)\\
&=\int_{\R} |x_1| \D\mu_1(x_1)+\int_{\R} |x_2| \D\mu_2(x_2)\\
&=\int_{\R^2}|x_1|+|x_2|\D\Q(x_1,x_2),
\end{align*}
which yields together with the weak convergence $\Q^{(n_k)} \rightarrow \Q$, by \cite[Definition 6.8.~(i)]{villani2008optimal}, that $\lim_{k \rightarrow\infty}\mathcal{W}\left(\Q^{(n_k)},\Q\right)=0$.
Next, let $H\in C_b(\R)$. We have that $(x_1,x_2) \mapsto H(x_1)(x_2-x_1) \in C_{\operatorname{lin}}(\R^2)$ and thus the Wasserstein convergence $\lim_{k \rightarrow\infty}\mathcal{W}\left(\Q^{(n_k)},\Q\right)=0$ implies, according to \eqref{eq_wassersteinspace_convergence}, that
\begin{align*}
0=\lim_{k\rightarrow\infty} \int_{\R^2} H(x_1)(x_2-x_1) \D \Q^{(n_k)}(x_1,x_2)=\int_{\R^2} H(x_1)(x_2-x_1) \D \Q(x_1,x_2).
\end{align*}
Since the function $H$ was chosen arbitrarily, we obtain by the definition of $\mathcal{M}(\mu_1,\mu_2)$, stated in \eqref{eq_martingale_def}, that $\Q \in \mathcal{M}(\mu_1,\mu_2)=\varphi\left((\mu_1,\mu_2)\right)$.
\item[3.)]
To show the lower hemicontinuity of $\varphi$, we apply Lemma~\ref{lem_lower_hemi}. Let $(\mu_1,\mu_2) \in X$ and consider a sequence $\left((\mu_1^{(n)},\mu_2^{(n)})\right)_{n \in \N} \subseteq  X$ such that 
\[
\lim_{n\rightarrow\infty}\mathcal{W}^{\oplus}\left((\mu_1^{(n)},\mu_2^{(n)}),(\mu_1,\mu_2)\right) =0.
\]
Note that, by step 1.), $\mathcal{M}(\mu_1,\mu_2)$ is non-empty.
Let $\Q \in \mathcal{M}(\mu_1,\mu_2)$. Then, by Theorem~\ref{thm_beigl_continuity}, there exists some sequence $(\Q^{(n)})_{n\in \N}$ with $\Q^{(n)} \in \mathcal{M}(\mu_1^{(n)},\mu_2^{(n)})$ for all $n \in \N$, converging w.r.t.\ $\mathcal{AW}$ towards $\Q$. Moreover, we have by definition the pointwise inequality $\mathcal{W} \leq \mathcal{AW}$. Thus, convergence w.r.t.\ $\mathcal{AW}$ implies convergence w.r.t.\ $\mathcal{W}$, i.e., we obtain
\[
\lim_{n\rightarrow\infty}\mathcal{W}\left(\Q^{(n)},\Q \right)\leq \lim_{n\rightarrow\infty}\mathcal{AW}\left(\Q^{(n)},\Q \right)= 0.
\]
\item[4.)] Let $(\mu_1^{(n)},\mu_2^{(n)},\Q^{(n)})_{n \in \N}\subseteq  \operatorname{Gr}\varphi$ be a sequence such that \[
\lim_{n\rightarrow\infty}\mathcal{W}^{\oplus}\left((\mu_1^{(n)},\mu_2^{(n)}),(\mu_1,\mu_2)\right)+\lim_{n\rightarrow\infty}\mathcal{W}\left(\Q^{(n)},\Q\right) = 0,
\]for some $(\mu_1,\mu_2,\Q) \in \operatorname{Gr}\varphi$.
Then, we  obtain by \eqref{eq_wassersteinspace_convergence} and since $\Phi \in C_{\operatorname{lin}}(\R^2)$ that 
\begin{align*}
\lim_{n\rightarrow \infty} f\left((\mu_1^{(n)},\mu_2^{(n)},\Q^{(n)})\right)&=\lim_{n \rightarrow \infty}\int_{\R^2} \Phi(x_1,x_2) \D \Q^{(n)}(x_1,x_2)\\
&=\int_{\R^2} \Phi(x_1,x_2) \D \Q(x_1,x_2)= f\left((\mu_1,\mu_2,\Q)\right).
\end{align*}
\end{itemize}
\end{proof}
\section*{Acknowledgements}
\noindent
Financial support by the Nanyang Assistant Professorship Grant (NAP Grant) \emph{Machine Learning based Algorithms in Finance and Insurance} is gratefully acknowledged. Both authors are very grateful for the helpful comments from two anonymous referees.
%
\bibliographystyle{plain}
\bibliography{literature}

\end{document}